\documentclass[11pt]{article}

\usepackage[all]{xy}

\newtheorem{tm}{Theorem}[section]
\newtheorem{lm}[tm]{Lemma}
\newtheorem{pr}[tm]{Proposition}
\newtheorem{rmk}[tm]{Remark}
\newtheorem{cor}[tm]{Corollary}
\newtheorem{ex}[tm]{Example}

\newtheorem{??}[tm]{Question}

\font\tenmsb=msbm10
\font\sevenmsb=msbm7
\font\fivemsb=msbm5
    
\newfam\msbfam
\textfont\msbfam=\tenmsb
\scriptfont\msbfam=\sevenmsb
\scriptscriptfont\msbfam=\fivemsb
\def\Bbb#1{{\fam\msbfam #1}}

\font\teneufm=eufm10
\font\seveneufm=eufm7
\font\fiveeufm=eufm5
\newfam\eufmfam
\textfont\eufmfam=\teneufm
\scriptfont\eufmfam=\seveneufm
\scriptscriptfont\eufmfam=\fiveeufm
\def\frak#1{{\fam\eufmfam\relax#1}}

\def\lorw{\longrightarrow}
\newcommand\n{\noindent}

\newcommand\ci{\cite}

\newcommand\rat{{\Bbb Q}}
\newcommand\comp{{\Bbb C}}

\newcommand\zed{{\Bbb Z}}

\newcommand\pn[1]{{\Bbb P}^{#1}}

\newcommand\blacksquare{{\hspace*{\fill} $\fbox{}$}}

\newcommand{\ke}{ \hbox{\rm Ker} }

\newcommand{\td}[1]{ \tau_{ \leq {#1} } }

\newcommand{\pe}{ {\cal P }  }

\title{The standard filtration on cohomology \\
with compact
supports \\
with an appendix on the base change map and 
\\
the Lefschetz
hyperplane theorem}
\author{
Mark Andrea A.  de Cataldo
\\ \, \\ \, \\ \,
{\em Dedicated to Andrew J. Sommese, on his 60th birthday} \\
{\em with admiration and respect}
}
\date{}

\begin{document}\maketitle

\begin{abstract}
We describe the standard and Leray filtrations
on the cohomology groups  with compact supports of a quasi projective variety
with coefficients in a constructible complex using flags of 
hyperplane sections on a partial compactification of 
a  related variety.
One of the key ingredients of the proof is the Lefschetz hyperplane theorem for perverse sheaves and,
in an appendix,
we discuss the  base change maps for 
constructible sheaves on algebraic varieties and
their role in a proof, due to Beilinson, of the Lefschetz hyperplane theorem.
\end{abstract}

\tableofcontents

\section{Introduction}
\label{intr}
Let $f: X \to Y$ be a map of algebraic varieties. The Leray filtration on  the 
(hyper)cohomology groups 
$H(X,\zed)= H(Y, Rf_*\zed_X) $ is  defined to be the standard filtration
on $H(Y, Rf_* \zed_X)$, i.e. the one  given by  the images in cohomology of the truncation maps
$\td{i} Rf_* \zed \to  Rf_* \zed$. Similarly, for  the 
cohomology groups  with compact supports $H_c(X,\zed)= H_c(Y, Rf_! \zed_X)$.

D. Arapura's paper  \ci{arapura} contains a geometric description of the Leray filtration on the cohomology groups  $H(X,\zed)$
for a proper map of quasi projective varieties $f:X \to Y$.
For example, if $Y$ is affine, then the
Leray filtration  is given, up to a suitable re-numbering,
 by the   kernels of the restriction maps
$H(X, \zed) \to H(X_i, \zed) $ to  a suitable collection of 
subvarieties  $X_i \subseteq X$.  This description implies at once
that the Leray filtration, in fact the whole Leray spectral sequence,
is in the category of mixed Hodge structures.

The same  proof works if we replace the sheaf  $\zed_X$
with  any  bounded  complex $C$  of sheaves of abelian groups on $X$
with constructible cohomology sheaves. Such complexes 
are simply called constructible.

Since  
the key constructions take place on $Y$,  given a constructible complex $K$ on $Y$,
one obtains  an analogous   geometric description
 for the standard filtration on $H(Y, K)$.
 For example, if  $Y$ is affine, then
 there is a collection of  subvarieties  $Y_i\subseteq Y$, obtained
as complete intersections
of  suitably  high degree hypersurfaces in special position,
such that the standard filtration is given by the
kernels of the restriction maps
 $H(Y, K) \to H(Y_i, K_{|Y_i})$. 

 The case of the Leray filtration for a proper map mentioned above
is then the  special case $K= Rf_* C$, and the varieties $X_i= f^{-1}(Y_i)$.  The properness
of the map is used to ensure, via  the proper base change theorem,
 that the natural base change maps
are isomorphisms,  so that, in view of the fact that
$H(X,C) = H(Y, Rf_*C)$,  we can identify the two maps
\[H(X,C) \lorw H(X_i, C_{|X_i}), \qquad H(Y,Rf_* C) \lorw
H(Y_i, Rf_*C_{|Y_i}),\]
and hence their kernels.

We do not know
of an analogous description of the Leray filtration
on the cohomology groups $H(X, C)$ for non proper maps $f: X \to Y$.

In \ci{arapura}, D. Arapura  also gives a geometric  description
 of  the Leray filtration  on the cohomology groups
 with compact supports
 $H_c(X,C)$ for a {\em proper} map $f: X\to Y$ of quasi projective varieties
 by first  ``embedding" the given morphism
 into  a morphism $\overline{f}: \overline{X} \to \overline{Y}$ of projective varieties,
 by identifying cohomology groups with compact supports
 on $Y$ with cohomology groups  on $\overline{Y}$,
 and then by applying his aforementioned result
 for   cohomology groups and proper maps.  In his approach, it is important that $f$ is proper,
 and the identity $f_! =f_*$ is used in an essential way.
 
 The purpose of this paper is to prove that,
 given a quasi projective variety $Y$ and a constructible complex $K$
 on $Y$ and, given a (not necessarily proper) map $f: X\to Y$ of algebraic varieties
 and a constructible complex $C$ on $X$,  one 
obtains a geometric description of the 
 standard filtration on the cohomology groups with compact supports
 $H_c(Y, K)$ (Theorem \ref{tmalce}), and of the  
 Leray filtration on the cohomology groups with compact supports
 $H_c(X, C)$ (Theorem \ref{tmalceste}).
  
The proof  still relies on the geometric description
 of  the cohomology groups $H(X,C)$ for   proper maps $f:X \to Y$.
 In fact, we utilize
 a completion $\overline{f}: \overline{X} \to \overline{Y}$
of  the varieties {\em and} of the map; see 
diagram (\ref{bd}). 

For completeness, we include
a new proof of the main result of \ci{arapura},
i.e. of  the geometric description of the Leray filtration on 
the cohomology groups $H(X,C)$ for proper maps $f:X \to Y$; see Corollary\ref{cortm1}.
In fact, we point out that we can extend the result to cover
the case of the standard filtration on the cohomology groups
$H(Y,K)$; see Theorem \ref{tm1}. Theorem \ref{tm1} implies
Corollary \ref{cortm1}.

The  proof of Theorem \ref{tm1} is based on the techniques
introduced in  \ci{decmigs1}, which deals
with perverse filtrations.  
In the perverse  case there is no formal difference
in the treatment of cohomology and  of cohomology with compact supports.
This contrasts sharply  with the standard case. 

Even though the  methods in this paper and in \ci{decmigs1, decII} are quite different
from the ones in \ci{arapura},  
the idea of  describing  filtrations geometrically by using hyperplane
sections comes from
 \ci{arapura}. 

In either approach,  the Lefschetz hyperplane
Theorem \ref{swls} for constructible sheaves on varieties with arbitrary singularities
plays a central role.
This result is due to several authors, Beilinson
\ci{be}, Deligne (unpublished) and Goresky and MacPherson \ci{gomasmt}.
Beilinson's proof works in the \'etale context and is a beautiful application
of the generic base change theorem.

In the Appendix \S\ref{bclht},  I discuss the  base change maps for 
constructible sheaves on algebraic varieties and the role played
by them in Beilinson's proof of the Lefschetz hyperplane
theorem. This is merely  an  attempt to make these
techniques more accessible to non-experts and hopefully
justifies 
the length of this section and the fact that it contains fact well-known
to experts.

 The notation employed in this paper is explained in
some detail, especially for non-experts, in 
$\S$\ref{notandback}. Here is a summary.
A variety is a separated scheme of finite type over $\comp$.
We work with 
 bounded complex  of sheaves of abelian groups on $Y$
 with constructible cohomology sheaves (constructible complexes, for short)
 and denote the corresponding derived-type category by ${\cal D}_Y$. 
The results hold, with essentially the same proofs, in the context of   \'etale cohomology
for varieties over algebraically closed fields; we do not discuss this variant. 
 For $K \in {\cal D}_Y$,   we have the (hyper)cohomology groups
 $H(Y, K)$ and $H_c(Y,K)$, the 
truncated complexes  $\tau_{\leq i} K$ and the cohomology sheaves
${\cal H}^i (K)$ which fit   into the exact sequences (or distinguished triangles)
\[ 0 \lorw \td{i-1} K \lorw \td{i} K \lorw  {\cal H}^i(K) [-i] \lorw 0. \]
Filtrations on abelian groups, complexes, etc., are taken to be decreasing,
$F^iK \supseteq F^{i+1}K$. The quotients (graded pieces)
are denoted $Gr^i_F K: = F^i K /F^{i+1}K$.

The standard (or Grothendieck) filtration on $K$ 
is defined by  setting $\tau^p K : = \td{-p}K$. The graded complexes
 satisfy $Gr^p_\tau K = 
{\cal H}^{-p}(K) [p]$. 
The corresponding decreasing and finite 
filtration  $\tau$ on the  cohomology   groups
$H(Y,K)$ and 
$H_c(Y,K)$ are 
called the standard (or Grothendieck) filtrations. 
Given a map $f: X \to Y$ and a complex $C \in {\cal D}_X$, the derived image complex
$Rf_* C$ and the derived image with proper supports complex $Rf_! C$
are in ${\cal D}_Y$ and 
the standard filtrations on
$H(Y, Rf_*C)=H(X, C)$ and $ H_c(Y, Rf_!C)=H_c(X,C)$ are called
the Leray filtrations.

A word of caution.
A key fact used in \ci{decmigs1} in  the case of the perverse filtration is that exceptional
restriction  functors $i^!_p$
to general linear sections $i_p:Y_p \to Y$ preserve perversity (up to a shift).
This fails in the case of the standard filtration
where we must work with   hypersurfaces in  special position.
In particular,
 $i^!$ of a sheaf is not a sheaf, even after a suitable shift, and this
 prohibits  the  extension of  our inductive approach in Theorem \ref{tm1}
 and in Corollary \ref{cortm1}  from cohomology
 to cohomology with compact supports.  As a consequence,
  the statements we prove
in cohomology for the standard and  and for the Leray filtrations  do not have a 
direct  counterpart in cohomology with compact supports, by, say,
a reversal of  the arrows. The remedy  to this offered in this paper
passe through  completions of
varieties and maps.

All the results of this paper  are stated in terms of filtrations on 
cohomology groups and on  cohomology groups with compact supports, 
but hold more generally,
and with the same proofs, 
for the associated filtered complexes and spectral sequences. However, for simplicity
of exposition,
we only state  and prove these results for  filtrations.

\medskip
{\bf Acknowledgments.} I thank Luca Migliorini for many conversations.

\section{The geometry of the standard and Leray filtrations}
\label{gslf}
In this section we give a geometric  description of standard and Leray filtrations
on cohomology and on cohomology with compact supports
in terms of flags of subvarieties.  
\subsection{Adapting  \ci{decmigs1} to the standard filtration
in cohomology: affine base}
\label{summ?}
While   the paper \ci{decmigs1} is concerned with the perverse filtration,
its  formal set-up  is quite general and is easily adapted to 
the case of the  standard filtration.  In this section, we briefly
go through the main steps of this adaptation and prove the key 
Theorem  \ref{tm1} and its Corollary \ref{cortm1}.
We refer the reader  to \ci{decmigs1} for more details and for the proofs
of  the various statements we discuss and/or list without proof.

The  shifted filtration $Dec(F)$   associated with a filtered complex of abelian 
groups
$(L,F)$ is  the filtration on $L$ defined as follows:  
\[ Dec(F)^n L ^l := \{ x \in F^{n+l}L ^l \, | \; dx \in F^{n+l+1}L^{l+1} \}.\]
The resulting filtrations in cohomology satisfy
\[ Dec(F)^n H^l(L) = F^{n+l} H^l(L).\]

\begin{pr}
\label{pr-1}
Let $(L, P, F)$ be a bifiltered complex of abelian groups. Assume that 
\begin{equation}
\label{-1}
H^r( Gr_P^b Gr_F^a L) =0  \quad
\forall \, r \neq a-b.
\end{equation}
Then
$L= Dec(F)$ on $H(L)$.
\end{pr}

Let $K \in {\cal D}_Y$ be a constructible complex on a variety $Y$.
 By replacing $K$ with a suitable injective resolution,
we may assume that $K$ is endowed with a filtration $\tau$ such that
the complex 
$Gr_\tau^b K [-b]$ is an injective resolution of the sheaf ${\cal H}^{-b}(K)$.
We take global sections and obtain the filtered complex
$(R\Gamma (Y, K), \tau)$ for which we have
$Gr^b_\tau R\Gamma (Y,K)= R\Gamma (Y, {\cal H}^{-b}(K)[b])$.
The filtration $\tau$ on $R\Gamma (Y,K)$ induces the standard filtration
(denoted again by $\tau$) on the
cohomology groups $H^*(Y,K)$.

An $n$-flag on $Y$ is an increasing sequence of closed subspaces
\[Y_*: \qquad  \emptyset = Y_{-1} \subseteq Y_{0} \subseteq \ldots \subseteq Y_n =Y.\]
The flag  $Y_*$   induces a filtration $F_{Y_*}$  on $K$ as follows:
(recall that $j_!= Rj_!$ and $k_!= Rk_!$ are extension by zero)
set  $j_a: Y\setminus Y_{a-1} \to Y$ 
and define  
\[ F^a_{Y_*} K := {j_a}_! j_a^* K= K_{Y-Y_a}. \]
Setting $k_a: Y_a\setminus Y_{a-1} \to Y$, we have 
$Gr_F^a K=  {k_a}_! k_a^* K= K_{Y_a -Y_{a-1}}$. 
The corresponding filtration in cohomology is 
\[
F^a_{Y_*} H^r (Y, K)  = \ke{\, \{H^r (Y, K) \lorw H^r (Y_{a-1}, K_{|Y_{a-1}}) \}}.\]

Taking global sections, we get the filtered complex
$(R \Gamma (Y, K), F_{Y_*})$ with the property
that 
\begin{equation}
\label{-2}
Gr_{F_{Y_*}}^a R\Gamma (Y, K) = R\Gamma (Y, K_{Y_a -Y_{a-1}})=
R\Gamma (Y_a, (K_{|Y_a})_{Y_a- Y_{a-1}}).
\end{equation}
We have 
\[
Gr^b_\tau Gr^a_{F_{Y_*}} R\Gamma (Y, K) =R\Gamma (Y_a, {\cal H}^{-b}(K) [b]_{Y_a -Y_{a-1}}),
\]
so that
\begin{equation}
\label{-111}
H^r (Gr^b_\tau Gr^a_{F_{Y_*}} R\Gamma (Y, K) )
=H^{r+b} (Y_a, ({\cal H}^{-b}(K)_{| Y_a}) _{Y_a -Y_{a-1}}).
\end{equation}

Note that the left-hand-side is
the relative cohomology group  
\begin{equation}
\label{relcg}H^{r+b} ( Y_a,  Y_{a-1},  {\cal H}^{-b}(K)_{|Y_a})=
H^{r+b}(Y_a,  {j_a}_! j_a^*{\cal H}^{-b}(K)_{|Y_a}),
\end{equation}
where $j_a: Y_a \setminus Y_{a-1} \to Y_a$. This is important in what follows
as it points to the use we now make of the Lefschetz hyperplane theorem for Sheaves \ref{swls}.

\begin{tm}
\label{tm1}
Let $Y$ be  an affine variety  of dimension $n$ and $K \in {\cal D}_Y$
be a constructible complex on $Y$. There is an $n$-flag $Y_*$ on $Y$
such that 
\[
\tau = Dec (F_{Y_*}) \qquad  \mbox{on $H(Y,K)$}.\]
\end{tm}
{\em Proof.}
The goal is to choose the flag  $Y_*$ so that (\ref{-1}) holds 
when $(L,P, F): = (R\Gamma (Y, K), \tau, F_{Y_*})$.
In view of (\ref{relcg}), we need the flag to satisfy the condition
\begin{equation}
\label{elccm}
H^r (Y_a, {j_a}_! j_a^* ({\cal H}^\beta (K))_{|Y_a}) =0 \qquad \forall \, r \neq a, \;\; \forall a \in [0,n],
\;\;  \forall \, \beta.
\end{equation}
Note that Theorem \ref{swls} applies to any finite collection of sheaves
(in fact it applies to any collection of sheaves which are constructible with respect to a fixed stratification). The flag is constructed by descending induction
on the dimension of $Y$. By definition, $Y_n =Y$. It is sufficient to choose
$Y_{n-1}$ as in Theorem \ref{swls}. We repeat this process, replacing $Y_n$ with $Y_{n-1}$
and construct the  wanted flag inductively.
\blacksquare

\medskip
Let $f:X \to Y$ be a map of algebraic varieties with $Y$ affine  and $C \in {\cal D}_X$. 
The Leray filtration $L^f$ on $H(X,C)=H(Y, Rf_* C)$ is, by definition,
the standard filtration $\tau$  on $H(Y,Rf_*C)$. 
Theorem \ref{tm1}  yields an $n$-flag $Y_*$ on $Y$ such that
$L^f= Dec(F_{Y_*})$.  In the applications though,
it is more useful to have a description in terms of a flag on $X$.
Let $X_* := f^{-1} Y_*$ be the pull-back flag on $X$, i.e. $X_a:= f^{-1} Y_a$.
There is the commutative diagram
\begin{equation}
\label{ttt}
\xymatrix{
H( {X} , \pi^* C) \ar[d]^{r}      & = &  &  H(   {Y}, Rf_* 
\pi^* C) \ar[d]^{r'} \\
H( {X}_{a}, i_a^* \pi^* C  )  & = &  H({Y}_a, Rf_* i_a^* \pi^* C) & 
H ( {Y}_a, i_a^* Rf_* \pi^* C) \ar[l]_{b}  ,
}
\end{equation}
where $b$ stems from the base change map (\ref{bc1})
$i_a^* Rf_* C \to Rf_* i_a^*C.$
The kernels
of the vertical restriction  maps $r$
and $r'$ define the filtrations $F_{X_*}$ and $F_{Y_*}$.
It is clear that $\ke{\, r} \supseteq \ke{\, r'}$, i.e. that
$F_{X_*}\supseteq F_{Y_*}$, and that  equality holds if
the base change map is an isomorphism. 

The following is now immediate.

\begin{cor}
\label{cortm1}
Let $f: X \to Y$ be a proper map with $Y$ affine of dimension $n$
and $C\in {\cal D}_X$. There is an $n$-flag $X_*$ on $X$ such that 
\[
L^f = Dec (F_{X_*}) \qquad  on \;\; H(X,C).
\]
\end{cor}
{\em Proof.} Since $f$ is proper, the base change map
$i_a^* Rf_* C \to Rf_* i_a^*C$ is an isomorphism.
\blacksquare

\medskip
In this section we have proved results for when $Y$ is affine.
In this case the statements and proofs are more transparent and
the flags are on $Y$ (pulled-back from $Y$ in the Leray case).
The case when $Y$  quasi projective case is 
easily reduced to the affine case in the next section.

\subsection{Standard and  Leray filtrations in cohomology: quasi projective
base}
\label{salfic}
In this section we  extend the results of the previous section
from the case when $Y$ is affine, to the case when $Y$ is quasi projective.
The only difference is that, given a quasi projective variety $Y$,
 we need to work with an auxiliary affine variety
${\cal Y}$ which is a fiber bundle, $\pi: {\cal Y} \to Y$, over $Y$ with fibers 
affine spaces ${\Bbb A}^d$ and we  need the flag to be an $(n+d)$-flag
${\cal Y}_*$ on ${\cal Y}$. This construction is due to Jouanolou.

Here is one way to prove this.  In the  case $Y = \pn{n}$,  take
${\cal Y} := (\pn{n} \times \pn{n}) \setminus \Delta$ with $\pi$
either projection. In general, take a projective completion $Y'$ of $Y$. 
Blow up  the boundary $Y' \setminus Y$ and obtain a projective completion
$\overline{Y}$ of $Y$ such that $Y \to \overline{Y}$ is affine.
Embed $\overline{Y}$ in some $\pn{N}$. Take the restriction
of the bundle projection $(\pn{N} \times \pn{N}) \setminus \Delta \to \pn{N}$
over $Y$ to obtain the desired result. 

Let $Y$ be a quasi projective variety. We fix  a ``Jouanolou  fibration"
$\pi: {\cal Y} \to Y$ as above.  If $Y$ is affine,  then we choose the identity.
In order to distinguish  standard filtrations
 on different spaces, e.g. $Y$ and ${\cal Y}$, we 
occasionally write
$\tau_Y$ and $\tau_{\cal Y}$.
Since the fibers of $\pi$ are contractible, we have canonical identifications of filtered groups
\[ (H(Y, K), \tau_Y)  =  (H({\cal Y}, \pi^*K), \tau_{\cal Y}). \]
This identity, coupled with Theorem \ref{tm1} yields
at once the following

\begin{tm}
\label{tm111}
Let $Y$ be quasi projective of dimension $n$ 
and $K \in {\cal D}_Y$. There is an $(n+d)$-flag ${\cal Y}_*$ on ${\cal Y}$
such that
\[ (H^*(Y, K), \tau_Y) = (H^*({\cal Y}, K), Dec(F_{{\cal Y}_*})), \]
i.e.
\[
\tau^p_Y \, H^r(Y, K)= \ke \, \{ \pi_{p+r-1}^*: H^r(Y,K) \lorw H^r ({\cal Y}_{p+r -1}, \pi_{p+r-1}^* K) \}.\]
If $Y$ is affine, then ${\cal Y} = Y$ and the flag is an $n$-flag on $Y$.
\end{tm}

In order to generalize Corollary \ref{cortm1} about the Leray filtration,
we form the Cartesian diagram (where maps of the ``same" type
are denoted with the same symbol)
\[
\xymatrix{
{\cal X} \ar[r]^f \ar[d]^{\pi} & {\cal Y} \ar[d]^\pi \\
X \ar[r]^f & Y }\]
 and  there are   the identifications (the last one stems from Base change for Smooth Maps)
\[H({\cal X}, \pi^* C)= H(X, C) = H(Y, Rf_*C) = H({\cal Y}, \pi^* Rf_* C)=
H({\cal Y}, Rf_* \pi^* C))\]
and we have the identity of the corresponding filtrations
\[
L^{f: {\cal X}\to {\cal Y}} =  L^{f: {X}\to {Y}} =\tau_Y = \tau_{\cal Y} = \tau_{\cal Y}.\]

\begin{tm}
\label{tm222}
Let $f: X \to Y$ be a proper map of algebraic varieties,  let $Y$ be quasi projective
 and $C \in {\cal D}_X$. There is an $(n+d)$-flag
${\cal X}_*$ on ${\cal X}$ such that
\[
(H(X,C), L^f) = (H({\cal X}, \pi^* C), Dec(F_{{\cal X}_*})),\]
i.e.
\[ L^p H^r(X,C) = \ke \,\{ \pi^*_{p+r-1}: H^r(X,C) \lorw H^r({\cal X}_{p+r-1}, 
\pi^*_{p+r-1} C)\}. \]
\end{tm}
{\em Proof.} We have $L^{f: {X}\to {Y}} = L^{f: {\cal X}\to {\cal Y}}$ so that
we may replace $f: X\to Y$ with $f: {\cal X} \to  {\cal Y}$. Since
${\cal Y}$ is affine,  we can apply Corollary \ref{cortm1} and conclude.
\blacksquare

\begin{rmk}
\label{rmknotpro}
{\rm 
If the map $f$ is not proper, then the relevant base change map is not an isomorphism.
While it is possible to describe the Leray filtration using a flag on ${\cal Y}$ (on $Y$ if $Y$ is affine), we do not know how to describe it using a
flag on ${\cal X}$ (on $X$ if $Y$ is affine). This latter description
would be more desirable in view, for example, of  the following
Hodge-theoretic application due to Arapura \ci{arapura}.
We also  do not know how to do so using compactifications; see 
Remark \ref{nonfunza}.
}
\end{rmk}

\begin{cor}
\label{corar}
Let $f:X \to Y$ be a proper morphism with $Y$ quasi projective.
Then the Leray filtration on $H(X,\zed)$ is by mixed Hodge substructures.
\end{cor}
{\em Proof.}  The Leray filtration satisfies $L^f=Dec( F_{{\cal Y}_*})$ for some 
flag on the auxiliary space ${\cal Y}$. Since the base change maps are isomorphisms,
$F_{ {\cal Y}_*} = F_{ {\cal X}_*}$. By the usual functoriality property
of  the canonical mixed Hodge structures on varieties, 
the latter filtration is given by  mixed Hodge substructures of 
$H({\cal X}, \zed) = H(X, \zed)$.
\blacksquare

\subsection{Standard and Leray filtrations in cohomology
with compact supports}
\label{salfccs}
Since  sheaves do not behave well with respect to Verdier Duality,
it is not possible to dualize the results in cohomology to obtain analogous ones for
cohomology with compact supports.

Given a map  of varieties $f:X \to Y$ and $C \in {\cal D}_X$, the Leray
filtration $L^{f:X \to Y}$ on $H_c(X, C)= H_c(Y, Rf_! C)$ is defined to be the standard filtration
on the last group. 

In this section we give a geometric description of the 
standard and of the Leray filtrations on cohomology with compact supports.
The description of the Leray filtration  on the cohomology groups  with compact supports $H_c(X,C)$ 
is valid for any (not necessarily proper)  map,
and this is in contrast with
the case of the  cohomology groups $H(X,C)$  (see Remark
\ref{rmknotpro}).

Arapura's \ci{arapura} proves an analogous result for proper maps, but to our knowledge,
that method does not extend to non proper maps. Nevertheless,
the method presented here is close in spirit to Arapura's.

The method consists of passing to completion of varieties and maps
and then 
use the 
 base change properties
associated with these compactifications to reduce to the case of cohomology and proper maps.
One main difference with cohomology is that, even if we start with $Y$ affine, the flag 
is  on the   auxiliary space $\overline{\cal Y}$ (see  below).
 
 \subsubsection{Completions of varieties and maps}
 \label{cvm}
We use freely the fact, due to Nagata, that varieties
and maps  can be compactified, i.e. any variety $Y$ admits an open immersion
into a complete variety with Zariski dense image, and given  any map
$f: X \to Y$, there are a  proper map $f' : {X'} \to Y$ and 
an open immersion $X \to X'$ with Zariski-dense image such that
$f'_{|X} = f$.

Since our result are valid without the quasi projectivity
assumption on $X$, we invoke Nagata's deep result.
If $X$ and $Y$ are both quasi  projective, then it is easy to get 
 by taking projective completions of $X$ and $Y$ and by  resolving the indeterminacies.

Let $f: X \to Y$ be a map with  $Y$  quasi projective.
Choose a projective completion $j : Y \to \overline{Y}$ such that
$j$ is an affine open embedding. This can be achieved by
first taking any projective completion and then by blowing up the
boundary. Choose an
${\Bbb A}^d$-fibration  $\pi : \overline{\cal Y} \to \overline{Y}$.
 Choose a completion $j: X \to \overline{X}$
such that $f$ extends to a (necessarily) proper $\overline{f}: \overline{X} \to \overline{Y}.$
Choose  closed embeddings $i: \overline{\cal Y}_a \to \overline{\cal Y},$
e.g. the constituents of a flag $\overline{\cal Y}_*$ 
 on $\overline{\cal Y}$.

There is the following commutative diagram, where the squares
and parallelograms labelled $\bigcirc \!\!\!\!\!1, \ldots, 
\bigcirc\!\!\!\!\!8\,$ are Cartesian:

\begin{equation}
\label{bd}
\xymatrix{
{\cal X}_a \ar[rrrrrrrr]^j \ar[dd]^i  \ar[ddddddrrr]_f
  &&&&&&&& {\overline{\cal X}}_{a} \ar[dd]_i 
\ar[ddddddlll]^{\overline{f}} 
\\
&&&& \bigcirc \!\!\!\!\!8  &&&&  \\
{\cal X} \ar[rrrrrrrr]^j \ar[dd]^\pi 
\ar[ddddddrrr]^{\stackrel{\mbox{$\bigcirc \!\!\!\!\!5 $}}\,}_f
  &&&&&&&& {\overline{\cal X}} \ar[dd]_\pi 
\ar[ddddddlll]_{\stackrel{ \mbox{ $\bigcirc \!\!\!\!\!6 $}  }{\,} }^{\overline{f}}
 \\
&&&& \bigcirc \!\!\!\!\!7  &&&& \\
X \ar[rrrrrrrr]^j  
\ar[ddddddrrr]^{\stackrel{\mbox{$\bigcirc \!\!\!\!\!3 $}}\,}_f 
&&&& &&&& 
\overline{X}
\ar[ddddddlll]_{\stackrel{ \mbox{ $\bigcirc \!\!\!\!\!4 $}  }{\,} }^{\overline{f}} 
\\
&&&& &&&& \\
&&& {\cal Y}_a \ar[rr]^j \ar[dd]^i   && {\overline{\cal Y}}_{a} \ar[dd]_i &&&  \\
&&&& \bigcirc \!\!\!\!\!2 &&&& \\
&&& {\cal Y} \ar[rr]^j \ar[dd]^\pi   && {\overline{\cal Y}} \ar[dd]_\pi &&&  \\
&&&& \bigcirc \!\!\!\!\!1 &&&&  \\
&&& Y \ar[rr]^j    && \overline{Y}. &&&
}
\end{equation}
If $f:X \to Y$
is not proper, then the  commutative trapezoids are not Cartesian.

Note that since $j: Y \to \overline{Y}$ and $\overline{\cal Y}$ are affine,
the map  $\pi: {\cal Y} \to Y$ is an ${\Bbb A}^d$-bundle with affine total space ${\cal Y}.$

We shall use freely the facts that follow.

\begin{enumerate}

\item  Due to the smoothness of the maps $\pi$ and
the properness of the maps $\overline{f},$
the base change Theorem holds for $\bigcirc \!\!\!\!\!1\,,$
$\bigcirc \!\!\!\!\!3\,,$
$\bigcirc \!\!\!\!\!4\,,$
$\bigcirc \!\!\!\!\!6$ and
$\bigcirc \!\!\!\!\!7\,.$

\item
For the remaining Cartesian squares $\bigcirc \!\!\!\!\!2\,,$
$\bigcirc \!\!\!\!\!8\,,$ and 
$\bigcirc \!\!\!\!\!5$ we have the base change maps:
\begin{equation}
\label{dbcm}
 i^* Rj_* \lorw R j_* i^*,  \qquad  i^*Rf_* \lorw Rf_* i^* .
\end{equation}

\item
The exactness of $j_!$ implies that $R^0j_! = j_!$ is simply extension by zero
and that 
it commutes with ordinary truncation,
i.e.  $j_! \circ  \td{i} = \td{i} \circ  j_!\,$; similarly, for the formation of cohomology
sheaves. 
The compactness of
$\overline{Y}$ and $\overline{X}$ implies that,
$H(\overline{Y}, -) = H_c(\overline{Y}, -)$, etc.  Recall that $H_c(Y, -) = H_c(\overline{Y},
j_!(-))$.
It follows that we have, for every $K \in  {\cal D}_Y,$ canonical identifications of filtered groups 
\begin{equation}
\label{beeq}
(H_c(Y, K), \tau_Y)= (H_c(\overline{Y}, j_!K), \tau_{\overline{Y}})=
(H(\overline{Y}, j_!K), \tau_{\overline{Y}}). \end{equation} 
If   $\overline{K}\in {\cal D}_{\overline{Y}}$ is 
{\em any} extension  of $K \in {\cal D}_Y$ to $\overline{Y}$, e.g. $j_! K$,  $j_*K$ etc.,
then 
we also have 
\begin{equation}
\label{beeq1}
(H_c(Y, K), \tau_Y)=
(H(\overline{Y}, j_! j^! \overline{K}), \tau_{\overline{Y}}).
\end{equation}
 Similarly, for the other open immersions $j$
in  diagram (\ref{bd}).
In view of 
the definition of relative cohomology
as the hypercohomology of $j_! j^! (-)$, 
we also have
 $H(\overline{Y}, j_! j^! \overline{K}) =(H(\overline{Y}, \overline{Y}\setminus Y; \overline{K})$.

\item
 There are canonical identifications: 
\begin{equation}
\label{sfok}
(H(Y, K), \tau_Y) = 
(H({\cal Y}, \pi^* K), \tau_{\cal Y}).
\end{equation}
 Similarly, for the other maps $\pi.$

\item There are canonical identifications: 
\begin{equation}
\label{sfok2}
(H(X, C), L^{f}_{\tau}) = (H(Y, f_* C), \tau_Y)= 
(H({\cal Y}, \pi^*  f_*C= f_* \pi^* C), \tau_{\cal Y})=
( H ( {\cal X}, \pi^* C ), L^{  f}_\tau).
\end{equation}
 Similarly, for the  $\overline{f}$ in $\bigcirc \!\!\!\!\!4$.

\item
Recall that for a filtration
$F$, the filtration $F(l)$ is defined by setting  $F(l)^i: = F^{l+i}$).
Since $\pi^* = \pi^! [-2d]$ is an exact functor, we have 
canonical identifications of filtered groups
\begin{equation}
\label{sfok3}
(H_c(X, C), L^f_{\tau}) = (H_c(Y, f_! C), \tau_Y)= 
(H_c({\cal Y}, \pi^!  f_!C = f_! \pi^! C), \tau_{\cal Y}(2d)).
\end{equation}
 Similarly, for the other maps of type $f$ and $\overline{f}.$

\end{enumerate}

\subsubsection{Filtrations on $H_c$ via compactifications}
\label{lll}
Let $Y$ be a  quasi projective variety  of dimension $n$ and $K \in {\cal D}_Y$
be a constructible complex on $Y$.
Consider any diagram  as in  (\ref{bd}).$\bigcirc \!\!\!\!\!1$.
We can choose $\overline{Y}$ to be of dimension $n$.   
\begin{tm}
\label{tmalce}
There is an  $(n+d)$-flag $\overline{\cal Y}_*$ on $\overline{\cal Y}$ for which we have the following identity of filtered groups
\[
(H_c(Y, K), \tau_Y) =H(\overline{\cal Y}, \pi^* j_! K), Dec (F_{\overline{\cal Y}_*})).
\]
\end{tm}
{\em Proof.} In view of   (\ref{beeq}) and of  (\ref{sfok}) applied to $\pi: 
\overline{\cal Y} \to \overline{Y}$,
we have canonical identifications of  filtered groups
\begin{equation}
\label{plmmlp}
(H_c(Y, K), \tau_Y) = (H ( \overline{Y}, j_! K), \tau_{\overline Y})=
(H( \overline{\cal Y}, \pi^* j_! K) =  H( \overline{\cal Y},   j_! \pi^* K), \, 
\tau_{\overline{\cal Y}}).
\end{equation}
The conclusion follows
from Theorem  \ref{tm111} applied to the pair 
$(\overline{Y}, j_! K)$.
\blacksquare

\bigskip
Let $f: X \to Y$ be a map of algebraic varieties and $C \in {\cal D}_X$.
Consider any diagram as in (\ref{bd}). We can choose $\overline{Y}$
to be of dimension $n$.
There are natural identifications of groups
\begin{equation}
\label{-ar}
H_c(X,C) =  H_c(\overline{X}, j_! C) = H(\overline{X}, j_!C) =
H(\overline{\cal X}, \pi^* j_!C =j_! \pi^* C).
\end{equation}
We  endow  these groups
with the  respective Leray filtrations. Note that since $\overline{f}$
is proper, the Leray filtration on $H_c(\overline{X}, j_!C)$ coincides
with the ones for $H(\overline{X}, j_! C)$.

\begin{tm}
\label{tmalceste}
There is an  $(n+d)$-flag $\overline{\cal X}_*$ on $\overline{\cal X}$ for which we have the following identity of filtered groups
\begin{equation}
\label{x6}
L^{f : X \to Y}_\tau  \;  = \;  L^{\overline{f} : \overline{X} \to \overline{Y}}_\tau \; = \;
Dec (F_{ \overline{\cal Y}_* }) \; = \; Dec(  F_{ \overline{\cal X}_* }),
\qquad \mbox{on} \; \;\;\;
H_c (X, C).
\end{equation}
\end{tm}
{\em Proof.} 
The filtration $L^{f:X \to Y}$ is the standard filtration
on $H_c(Y, Rf_! C)$ which in turn, by the exactness of $j_!$
and the equality $H_c(Y, -)= H_c (\overline{Y}, j_! (-))$,
 coincides with the standard
filtration on $H_c(\overline{Y}, j_! Rf_! C)$. Since $\overline{Y}$ is compact
and $\overline f$ is proper (so that $Rf_! = Rf_*$), by the commutativity
of the base trapezoid diagram in (\ref{bd}), we have that
$R\overline{f}_! j_! = j_! Rf_!$ so that
$H_c (\overline{Y}, j_! Rf_! C)=H (\overline{Y}, Rf_*j_! C)$.
This implies the equality $L^{f : X \to Y}_\tau  \;  = \;  L^{\overline{f} : \overline{X} \to \overline{Y}}_\tau$. We are now in the realm of cohomology and proper maps
and the rest follows from Theorem \ref{tm222}  applied to $\overline{f}$
and to $j_!C$.
\blacksquare

\begin{cor}
\label{ht}
The Leray filtration on $H_c(X, \zed)$ is by mixed Hodge substructures.
\end{cor}
{\em Proof.}  By Theorem \ref{tmalceste} and (\ref{beeq1}),
 the filtration in question is 
 the  one induced by the flag $\overline{\cal X}_*$ on the relative cohomology
group $H(\overline{\cal X}, j_!  \zed)= H( \overline{\cal X}, {\cal X}, \zed)$.
The result follows from Deligne's mixed Hodge Theory \ci{ho3}.
\blacksquare

\begin{rmk}
\label{rmkcorara}
{\rm
The case when $f$ is proper is proved in \ci{arapura}.
}
\end{rmk}

\begin{rmk}
\label{nonfunza}
{\rm
If one tries to imitate the procedure we have followed in the case of  
the cohomology groups with   compact supports $H_c(X,C)$ for an arbitrary map
$f: X \to Y$,
with the goal of obtaining an
analogous result for  the Leray filtration on the cohomology groups
$H(X,C)$, 
then one hits the following obstacle: indeed, 
there are identifications $H(Y,Rf_*C) = H(\overline{Y}, Rj_*  R\overline{f}_* C) = 
H(\overline{Y},   R\overline{f}_*   Rj_*  C)$, however, since
$Rj_*$ does not commute with truncation, the  Leray filtrations
for $f$ and $\overline{f}$
do not coincide, and the imitation of the procedure would
yield a geometric description only for the case of  $\overline{f}$.

}
\end{rmk}

\section{Appendix: Base change and Lefschetz hyperplane theorem}
\label{bclht}

\subsection{Notation and background results}
\label{notandback}
$\,$

{\bf Varieties and maps.} 
A {\em variety} is a separated scheme of finite type over the field
of complex numbers $\comp$. In particular, we do not assume that
varieties are irreducible, reduced, or even pure dimensional.
Since we work inductively with intersections
of special hypersurfaces, we need this generality even if we start with a nonsingular
irreducible variety.   A {\em map} is a map of varieties, i.e.
map of $\comp$-schemes. 

{\bf Coefficients.}
The results
of this paper hold for sheaves of $R$-modules,
where $R$ is a commutative ring with identity
with finite global dimension, e.g. $R= \zed,$ $R$ a field, etc.
For the sake of exposition
 we work with $R=\zed$, i.e. with sheaves of abelian groups.
 
 {\bf Variants.}
  The results of this paper hold, with routine
 adaptations of the proofs, in the case
 of varieties over an algebraically closed field
 and \'etale sheaves with the usual
 coefficients: $\zed/l^m\zed,$ $\zed_l,$ $\rat_l$,
 $\zed_l [E]$, $\rat_l [E]$ ($E \supseteq \rat_l$ a finite extension)
 and $\overline{\rat}_l.$
 These variants  are not discussed further (see  \ci{bbd}, $\S$2.2 and $\S$6).

{\bf Stratifications.} 
  The term {\em stratification} refers
 to  an  algebraic Whitney stratification \ci{borel, goma2, gomasmt}.
 Recall that any two stratifications admit a common refinement
 and that maps of varieties can be stratified.  See also $\S$\ref{susec-001}.
 
 {\bf The constructible derived category ${\cal D}_Y$.}
 Let $Y$ be a variety,  
 $Sh_Y$  be the abelian category of sheaves of abelian 
 groups on $Y$ and $D (Sh_Y)$ be  the associated derived category.
 A sheaf $F \in Sh_Y$ is {\em constructible} if there is a 
 stratification of $Y$ such that the restriction
 of $F$ to each stratum is locally constant with stalk a finitely generated
 abelian group. 
 A complex is {\em bounded} if the cohomology sheaves
 ${\cal H}^i(K)=0$ for $|i| \gg 0$. A complex $K \in D(Sh_Y)$ with   constructible
 cohomology sheaves is said to be {\em constructible}.
 The category
 ${\cal D}_Y= {\cal D}_Y(\zed)$ is the full subcategory of the 
  derived category  $D(Sh_Y)$ whose  objects are the bounded  
  constructible complexes. For a given stratification $\Sigma$ of $Y$, a complex
 with this property is called $\Sigma$-constructible.
 Given a stratification $\Sigma$ of $Y$, there is the
full subcategory
${\cal D}_Y^{\Sigma} \subseteq {\cal D}_Y$   of complexes which are 
$\Sigma$-constructible. Hypercohomology groups are denoted  $H(Y,K)$
and   $H_c(Y,K)$.
If $K \in {\cal D}_Y$ and $n \in\zed$, then $K[n] \in {\cal D}_Y$ is the
({\em $n$-shifted}) complex with $(K[n])^i= K^{i+n}$. One has, for example,
$H^i(Y,K[n]) = H^{i+n}(Y,K)$.

{\bf The four functors associated with a map}.
Given a map  
$f:X \to Y$,  there are the  usual four  functors $(f^*,  Rf_*,Rf_!, f^!)$.  By abuse of  notation,  denote $Rf_*$ and $Rf_!$
simply by $f_*$ and $f_!$.
The four functors
preserve stratifications, i.e. 
if $f: (X, \Sigma') \to (Y, \Sigma)$ is stratified, 
then  $f_*, f_!: D^{\Sigma'}_X \to D^{\Sigma}_Y$
and 
$f^*, f^!:  D^{\Sigma}_Y \to D^{\Sigma'}_X.$

{\bf Verdier Duality.} The Verdier Duality functor ${\Bbb D}= {\Bbb D}_Y : {\cal D}_Y \to {\cal D}_Y$
is an autoequivalence with ${\Bbb D}\circ {\Bbb D} = \mbox{Id}_{{\cal D}_Y}$ and it  preserves stratifications. We have
${\Bbb D}_Y f_! =  f_* {\Bbb D}_X$ and
${\Bbb D}_X f^! =  f^* {\Bbb D}_Y$.

 {\bf Perverse sheaves.}
 We consider only  the {\em middle} perversity $t$-structure
 on ${\cal D}_Y$ \ci{bbd}.  
 There is  the  full subcategory $\pe_Y \subseteq {\cal D}_Y$
 of perverse sheaves on $Y$. The elements are special
 complexes in ${\cal D}_Y$. An important example is the intersection complex
 of an irreducible variety \ci{goma2, borel}.  Let $j : U \to Y$ be an open immersion;
 then $j^!=j^*: \pe_Y \to \pe_U$, i.e. they preserve perverse sheaves.
 Let $j: U \to Y$ be an {\em affine} open immersion; 
 then $j_!, j_*:\pe_U \to \pe_Y$. 
The Verdier Duality functor ${\Bbb D}: \pe_Y \to \pe_Y$
 is an autoequivalence.

{\bf Distinguished triangles for a locally closed   embedding.}
There is the notion of distinguished triangle in ${\cal D}_Y$: it is a sequence
of maps $X \to Y \to Z \to X[1]$ which is isomorphic in ${\cal D}_Y$
to the analogous sequence of maps arising from the cone construction associated with
a map of complexes $X' \to Y'$. 
Let $j: U \to Y $ be a locally closed  embedding  with associated
``complementary"
 embedding and $i: Y \setminus U \to Y$.
For every $K \in {\cal D}_Y,$ we have  distinguished triangles
 \begin{equation}
\label{adt11}
j_!j^! K \lorw K \lorw i_* i^* K \stackrel{[1]}\lorw, \qquad
i_!i^! K \lorw K \lorw j_* j^* K \stackrel{[1]}\lorw.
\end{equation}

{\bf Various base change maps.}
Given  two maps $ Y' \stackrel{g}\to Y \stackrel{f}\leftarrow X,$
there is the Cartesian diagram 
\begin{equation}
\label{cdm}
\xymatrix{
X' \ar[r]^{g} \ar[d]^{f} & X \ar[d]^f \\
Y' \ar[r]^g & Y.
}
\end{equation}
The ambiguity of the notation (clearly the two maps $g$ are different
from each other, etc.) does not generate ambiguous
statements in what follows, and  it simplifies the notation.

\n
There are the natural maps
\begin{equation}
\label{nmnbc}
g_! f_* \lorw f_* g_!, \qquad g^* f^! \lorw f^! g^*.
\end{equation}
There are  the base change maps
\begin{equation}
\label{bc1}
g^* f_* \lorw  f_* {g}^*, \qquad g^! f_! \longleftarrow f_! g^!
\end{equation}
and the base change isomorphisms
\begin{equation}
\label{bc2}
g^* f_! \stackrel{\simeq}\lorw  f_! {g}^*, \qquad  g^! f_* \stackrel{\simeq}\longleftarrow  f_* g^!.
\end{equation}
Similarly, for the higher direct images $R^if_*$ and $R^if_!$.
\begin{ex}
\label{fiber}
{\rm 
Let $Y'\to Y$ be the closed embedding of a point $y \to Y$. 
The first  base change map in (\ref{bc1})   yields a map $ 
(R^if_* \zed_X)_y \to H^i (f^{-1}(y), \zed)$.  Given a sufficiently small, contractible
neighborhood of $y$ in $Y$, we have $H^i(f^{-1}(U_y),\zed) = 
(R^if_* \zed_X)_y$. This base change map is seldom an isomorphism, e.g.
the open immersion
$X =\comp^* \to \comp =Y$, $y =0$.
This  failure  
   is corrected  in  (\ref{bc2}) by taking  the direct image
with proper supports. }\end{ex}

{\bf Base change theorems.}
The base change maps (\ref{bc1}) are  isomorphisms if either
one of the following conditions is met:
$f$ is proper, 
$f$ is locally topologically trivial over $Y$,
or $g$ is smooth.

{\bf The octahedron axiom.}
This is one of the axioms for a triangulated category and
  can be found in \ci{bbd}, 1.1.6. Here is
a convenient way to display it (see \ci{bbd}, 1.1.7.1). Given a composition 
$X \stackrel{f}\to Y \stackrel{g}\to Z$ of morphisms
 one has the following diagram
\begin{equation}\label{oct}
\xymatrix{
&& Z' \ar[rd] && \\
& Y \ar[ru] \ar[r]^g &Z \ar[r] \ar[rrd] &Y' \ar[rd]& \\
X\ar[ru]^f \ar[rru]_{gf} &&&& X'
}
\end{equation}
where $(X,Y,Z'), (Y,Z, X'), (X,Z,Y')$ and $(Z',Y',X')$ are distinguished triangles.
\begin{rmk}
\label{u-1}
{\rm
It is clear that  $g$  is an isomorphism if and only if  $Y' \simeq 0$ if and only if
$X' \to Z'[1]$  is an isomorphism.
}
\end{rmk}

{\bf The term ``general."}
Let  ${\frak P}$ be a property  expressed in terms of the hyperplanes
of a projective space  ${\Bbb P}$, i.e. the elements of $\pn{\vee}$.
We say that   {\em property ${\frak P}$ holds for a   general hyperplane} if there is a Zariski-dense open
subset   $V \subseteq \pn{\vee}$ such that   property
${\frak P}$ holds for every hyperplane in $V$. Of course this terminology
applies to propositions ``parameterized" by irreducible varieties and 
one can talk about a general pair of hyperplanes, in which case the variety
is $\pn{\vee}\times \pn{\vee}$, etc.
 
\medskip
As it is customary, we often denote a canonical isomorphism
with the symbol ``$=$."

\subsection{Base change with respect to subvarieties}
\label{bcintro}
Let us discuss  the following  two special cases
of (\ref{cdm}). Even though (\ref{cdm-1}) is a special case
of (\ref{cdm-000}),  it is convenient to distinguish between the two
(see Propositions \ref{precise01}, \ref{precise02}). Let $i: H \to Y$
be a closed embedding. Let $j: U\to Y$ be an open embedding
and $f:X \to Y$ be a map. We obtain the following two Cartesian diagrams
\begin{equation}
\label{cdm-1}
\xymatrix{
U\cap H \ar[r]^{i} \ar[d]^{j} & U \ar[d]^j \\
H \ar[r]^i & Y,
}
\end{equation}
\begin{equation}
\label{cdm-000}
\xymatrix{
X_H\ar[r]^{i} \ar[d]^{f} & X \ar[d]^f \\
H\ar[r]^i & Y.
}
\end{equation}

\begin{??}
\label{q-0}
{\rm 
Let $K \in {\cal D}_Y$, $C \in {\cal D}_X$.
Which conditions on the closed embedding $i:H \to Y$ ensure
that   the base change maps
\[ j_* i^* K \lorw i^* j_* K, \qquad
 j_! i^! K \longleftarrow i^! j_! K \]
 are isomorphisms? 
 
\n
Which conditions on the closed embedding $i:X_H \to X$ ensure
that   the base change maps
\[
f_* i^* K \lorw i^* f_* K, \qquad
f_! i^! K \longleftarrow i^! f_! K\]
are isomorphisms?
}
\end{??}

Answers to  these questions  are given in Propositions 
\ref{precise01}, \ref{precise02}. Since these results
 involve the notion
of stratifications, we discuss  briefly stratifications in the next section.

\subsubsection{Some background on stratifications}
\label{susec-001}
$\,$

{\bf Stratifications.} For background, see
\ci{borel, gomasmt}.  The datum of a stratification $\Sigma$
of the variety $Z$ includes a disjoint union
decomposition $Z = \coprod \Sigma_i$ into locally closed
nonsingular {\em irreducible} subvarieties  $\Sigma_i$
called {\em strata}. One requires that the closure
of a   stratum is a union of strata. These data are subject to the Whitney
Conditions A and B, which we do not discuss  here.
Every variety admits a stratification. Any two stratifications of the same variety
admit a common refinement.  Given a stratification $\Sigma$ of $Z$,
every point $z \in Z$ admits a fundamental system of {\em standard neighborhoods}
homeomorphic, in a stratum-preserving-way, to 
 ${\comp}^l \times {\cal C}({\cal L})$, where ${\cal C}$ denotes the real cone
 (with vertex $v$),
 ${\cal L}$ is the {\em link} of $z$ in $Z$ (relative to $\Sigma$) (it is a stratified
 space obtained by embedding $Z$ in some manifold, intersecting
 $Z$ with a submanifold meeting the stratum transversally at   $z$ 
 and  then  intersecting the result  with a small ball centered at $z$), and ${\comp}^l\times v$
 is the intersection of the stratum $\Sigma_i$ to which $z$ belongs
 with a small ball (in the big manifold containing $Z$) centered at $z$.

{\bf Constructible complexes.} 
Let $\Sigma$ be a stratification of  the variety $Z$.
The bounded complexes constructible with respect to $\Sigma$
form the category $D^{\Sigma}_Z $ which is a full
subcategory of the constructible derived category ${\cal D}_Z$.
If $K \in D^{\Sigma}_Z$,
$z \in Z$,  $U:= {\comp}^l \times {\cal C}({\cal L})$ is a standard neighborhood
of $z$ with second projection $\pi$, then $K_{|U} \simeq \pi^* \pi_* 
K_{|U}$, i.e. $K$ is locally a pull-back from the cone over the link.

{\bf Stratified maps.}  Algebraic maps can be stratified:
given  a map  $f: X \to Y$,
there are stratifications $\Sigma_X$ for $X$  and $\Sigma_Y$ for $Y$ such that 
(i) for every stratum  $S$   on  $Y$, the space $f^{-1} S $
is a union of strata on  $X$ and (ii) 
for  $y \in S$ 
 there exists a neighborhood $U$ of $y$ in $S$, a stratified space $F$
 and a stratification-preserving homeomorphism $F\times U \simeq f^{-1} U$
 which transforms the projection onto $U$ into $f$. If $f$ is a closed embedding,
 then each  stratum in $X$ is the intersection
 of $X$ with a   stratum of $Y$ of the same dimension.
 If $f$ is an open immersion, using standard neighborhoods,
  the local model  at $y \in Y \setminus  X$
 is $f: {\comp}^l \times ({\cal C}({\cal L})  -  {\cal C}({\cal L}'))
 \to  {\comp}^l \times {\cal C}({\cal L})$, where ${\cal L}'$
 is the link  at $z$ of $Y\setminus X$.

{\bf Normally nonsingular inclusions.} 
 A closed embedding $i: H \to Y$ is {\em normally nonsingular} with respect to a stratification $\Sigma$ of $Y$ if
$H$ is obtained locally on $Y$  by the following procedure: embed $Y$ into a manifold
$M$
and $H$ is the intersection $H' \cap Y$, where $H'\subseteq M$
is a submanifold meeting transversally all the strata of $\Sigma$. See \ci{goma2, decmightam}.
Note that a normally nonsingular inclusion is locally of pure
codimension.
 If $Y$ is embedded into
some projective space, then by the Bertini Theorem a general hyperplane section
yields a normally nonsingular inclusion.  More generally,  the general element
of a finite dimensional  base-point-free linear system of $Y$ yields a normally nonsingular inclusion (\ci{jobert}).
Let $K \in D^{\Sigma}_Y$ and $i:H \to Y$ be  a normally nonsingular
inclusion of complex codimension $r$ with respect to $\Sigma$.
Then $i^* K = i^! K [-2r]$ (cf. \ci{decmightam}). 

\subsubsection{Sufficient condition
for the base change maps to be an  isomorphism}
\label{susec-11}
Let $j: U \to Y$ be an open embedding and $\Sigma$ be a stratification
of $Y$ such that, if $\Sigma_U$ is its trace on $U$, then the map $j: (U, \Sigma_U)
\to (Y, \Sigma)$  is stratified. Such a  stratification $\Sigma$ exists.
Consider the situation (\ref{cdm-1}).  
\begin{pr}
\label{precise01} 
Assume that $i:H \to  Y$ is normally nonsingular
with respect to ${\Sigma}$.
Then for every  $K \in D^{\Sigma_U}_U$ the base change maps
$$
i^* j_* K \lorw j_* i^* K, \qquad 
i^! j_! K \longleftarrow j_! i^! K
$$
are  isomorphisms.
\end{pr}
{\em Proof.}
Here are   two essentially equivalent proofs.
While the first one seems shorter,  it does  rely on
the formula $i^* = i^! [-2r]$, the second one is more direct.

\n
{\em $1^{st}$ proof.} The assumptions imply that
$i^* K = i^! K [-2r]$. The conclusion follows
from the base change isomorphisms (\ref{bc2}).

\n
{\em $2^{nd}$ proof.}
The complexes $j_! K, j_* K \in D^{\Sigma }_{Y }$.
The assertion is local. The local model for (\ref{cdm-1})
at a point $y \in H$ lying on a $l$-dimensional stratum
with links ${\cal L}$ for $y \in Y$ and ${\cal L}' \subseteq
{\cal L}$ for $y \in Y\setminus U$ is, denoting by ${\cal C}$ real cones
and by $r$ the codimension of $H$ in $Y$:
\begin{equation}
\label{summ}
\xymatrix{
 {\comp}^{l-r} \times ({\cal C}( {\cal L}) \setminus  {\cal C}({\cal L}')   )
\ar[rr]^i \ar[d]^j  &&
{\comp}^{l} \times  ({\cal C}( {\cal L}) \setminus  {\cal C}({\cal L}')   )
\ar[rr]^\pi \ar[d]^j  &&
\{ y \} \times 
 ({\cal C}( {\cal L}) \setminus  {\cal C}({\cal L}')   )  \ar[d]^j   \\
{\comp}^{l-r} \times {\cal C}( {\cal L}    )
\ar[rr]^i  &&
{\comp}^{l} \times {\cal C}( {\cal L}   )
\ar[rr]^\pi   &&
\{ y \} \times 
{\cal C}( {\cal L}    )  
}
\end{equation}
with $Id \simeq \pi^*  \pi_*= \pi^! \pi_! $ for 
$\Sigma$-constructible complexes.
One has, using the  base change Theorem for the smooth
map  $\pi \circ i $: 
$$ i^* j_* K \simeq i^* \pi^* \pi_*  j_* K = i^* \pi^* j_* \pi_* K =
j_* i^* \pi^* \pi_* K \simeq j_* i^* K.$$
This proves the first assertion. The second one is proved
in the same way:
$$ i^! j_! K \simeq i^! \pi^! \pi_!  j_! K = i^! \pi^! j_! \pi_! K =
j_! i^! \pi^! \pi_! K \simeq j_! i^! K.$$
Note that once the first assertion is proved, one can   prove the second one also
 as follows.
Given $K \in D^{\Sigma_U}_U$, we have that 
$K^{\vee} \in D^{\Sigma_U}_U$. 
We have proved that the first assertion holds
for every $K \in D^{\Sigma_U}_U$ so that it holds for $K^{\vee}$:
$i^* j_* K^{\vee} \simeq j_* i^* K^{\vee}$ and the second assertion 
follows by applying Verdier Duality to this isomorphism.
\blacksquare

\bigskip
Consider the Cartesian diagram (\ref{cdm-000}) and
let $\Sigma'$ be a stratification of $X$.
\begin{pr}\label{precise02}
Let $i:X_H \to X$ be a normally nonsingular inclusion
with respect to $\Sigma'$. Then  for every $C \in {\cal D}_X^{\Sigma'}$
the base change maps
\[
i^*  f_* C \lorw  f_* i^* C, \qquad
i^!  f_! C \longleftarrow   f_! i^! C\]
are  isomorphisms.
\end{pr}
{\em Proof.}
 Let $j: X \to \overline{X} \stackrel{\overline{f}}\to Y$
be a completion of the map $f,$
i.e. $j$ is an open immersion with dense image
and $\overline{f}$ is proper. Such a completion exists by
a fundamental result of Nagata.  
The Cartesian diagram (\ref{cdm-000}) can be completed 
to a   commutative diagram
with Cartesian squares
\begin{equation}
\label{ddiiaagg}
\xymatrix{
X_H \ar[rrr]^i \ar[rd]^j \ar[ddd]^f  &&& X  \ar[rd]^j \ar[ddd]^f & \\
 & \overline{X}_H \ar[rrr]^i \ar[ddl]^{\overline f} &&& \overline{X} 
\ar[ddl]^{\overline f} \\
&&&& \\
H \ar[rrr]^i  &&& Y &.
}
\end{equation}
By virtue of Lemma  \ref{precise01}, we have the base change isomorphism
$
i^* j_* C \simeq j_* i^* C
$
to which we apply the proper map $\overline{f} : \overline{X}_H \to H$:
$$
 \overline{f}_* i^* j_* C \simeq  \overline{f}_* j_* i^* C = f_* i^* C.
$$
The first assertion follows by applying  the Base Change for Proper Maps to the
first term:
\[ \overline{f}_* i^*  j_* C =  i^*  \overline{f}_*  j_* C =i^* f_* C.\]
The second one can be proved in either of  two ways
as in  the proof of Proposition \ref{precise01}.
\blacksquare

\begin{rmk}
\label{nosame}
{\rm
Note that $X_H$ can be normally included in $X$ with respect $\Sigma'$
while $H$ may fail to be so with respect to any  stratification of $Y$ for $f_*C$ and $f_!C$.
E.g. $X$ is nonsingular, $C =\zed_X$, but $Y$ and/or $f$  have singularities.
}
\end{rmk}
\subsection{The Lefschetz hyperplane theorem}
\label{lhtbe}
The classical Lefschetz hyperplane theorem  states that
if $Y$ is a projective manifold of dimension $n$ and $H$ is a smooth
hyperplane section relative to an embedding into projective space,
then the restriction map $H^i(Y,\zed) \to H^i(H, \zed)$ is an isomorphism
for $i < n-1$ and injective for $i =n-1$. 

In fact, one can prove  this (see \ci{milnor}, for example) by showing
that  the relative cohomology groups $H^i(Y, H, \zed) =0$ for $i \leq n-1$.

If $j: Y \setminus H \to Y$ is the open embedding, then $j^* = j^!$ and,
since $Y$ is compact, 
$H^i (Y, H, \zed) =
H(Y, j_! j^* \zed_Y)$. We can thus    reformulate the Lefschetz hyperplane theorem
in terms of the following vanishing statement
\[
H^i (Y, j_! j^! \zed[n]) = 0  \quad \forall \; i  <0.
\]
Note that $\zed[n]$ is a perverse sheaf on  the nonsingular $Y$.
Beilinson \ci{be}, Lemma 3.3,  has given a proof 
of this important result which is valid in the \'etale case and for every perverse sheaf
on a quasi projective variety $Y$. His proof is based
on the   natural map (\ref{cesem})  being an isomorphism.

In this section, we discuss Beilinson's proof, which boils down
to  an application of the
base change Proposition \ref{precise01}.

\subsubsection{The natural map $j_! J_* \lorw J_* j_!$}
\label{ss-1} Let $Y$ be  a quasi projective variety, $Y \subseteq \pn{N}$
be a fixed embedding in some projective space, 
$\overline{Y}\subseteq \pn{N}$  be the closure of $Y$, $\Lambda\subseteq 
\pn{N}$ be a hyperplane  and 
 $\overline{H} \subseteq \overline{Y}$ and $H \subseteq Y$ be the corresponding hyperplane sections.
 There is  the Cartesian  diagram
\begin{equation}
\label{su}
\xymatrix{
 H \ar[r]^i \ar[d]^J & Y \ar[d]^J & U  \ar[d]^J \ar[l]_j \\
\overline{H} \ar[r]^i & \overline{Y} & \overline{U}  \ar[l]_j.
}
\end{equation}

Let $K \in {\cal D}_Y$.
Consider the composition
\begin{equation}
\label{eq-100}
 J_* K \lorw i_* i^* J_* K \stackrel{\phi}\lorw i_* J_* i^* K (=J_* i_* i^* K).\end{equation}
The octahedron axiom  yields a distinguished triangle (the equality  stems
from  (\ref{bc2}))
\begin{equation}
\label{eq-101} j_! j^! J_* K[1]  (= j_! J_* j^! K [1]) \lorw  J_* j_!j^! K [1] 
\lorw \mbox{Cone} (\phi) \lorw, 
\end{equation}
where the first map arises by applying  (\ref{nmnbc}) to $j^!K$.

Similarly, we have the composition
\begin{equation}
\label{eq-100bis}
 J_! K \longleftarrow i_! i^! J_! K \stackrel{\varphi}\longleftarrow i_! J_! i^! K 
 (=J_! i_! i^! K) 
 \end{equation}
and the octahedron axiom,  yields a distinguished triangle
\begin{equation}
\label{eq-101bis}
\longleftarrow j_* j^* J_! K  (= j_* J_! j^* K ) \longleftarrow  J_! j_*j^* K 
\longleftarrow \mbox{Cone} (\varphi), 
\end{equation}
where the second  map arises by applying  (\ref{nmnbc}) to $j^*K$.

\begin{lm}
\label{pap}
The map  
\begin{equation}
\label{eq-10}
\label{cesem}
 j_! J_* j^! K  \lorw J_* j_! j^! K  \qquad 
 ( j_* J_! j^* K  \longleftarrow J_! j_* j^* K, \, resp.)
 \end{equation}
is an isomorphism 
if and only if the base change map $i^* J_* K \lorw  J_* i^* K$ 
($i^! J_! K \longleftarrow  J_! i^! K$, resp.)
is an isomorphism.
\end{lm}
{\em Proof.} By Remark \ref{u-1}, the map (\ref{cesem}) is an isomorphism if and only if
the map $\phi$ is an isomorphism. The conclusion follows from the 
fact that since $i$ is a closed embedding, $i_*$ is fully faithful.
The second assertion is proved using the same
construction, with the arrows reversed.
\blacksquare

\begin{rmk}
\label{no-es-i}
{\rm If in the set-up of  diagram (\ref{su}) the maps  $J : Y \to \overline{Y}$
and $i : \overline{H} \to {\overline{Y}}$  are  
arbitrary locally closed embedding of varieties, then
the proof of  Lemma \ref{pap} shows that {\em if the base change map
$i^* J_* K \to J_* i^* K$ is an isomorphism, then
the map (\ref{cesem}) is also an isomorphism}.
}
\end{rmk}

\subsubsection{The Lefschetz hyperplane theorem for perverse sheaves}
\label{tlht}
In this section,
 $Y$ is a {\em quasi projective variety} equipped with a fixed {\em affine} embedding
$Y \subseteq \pn{N}$ is some projective space.   Let us stress
that we shall consider
hyperplane sections with respect to this fixed affine embedding.

If $Y$ is affine, then {\em every} embedding into projective space is affine.
Not every embedding of a quasi projective variety is affine,
e.g. ${\Bbb A}^2 \setminus \{(0,0) \}\subseteq \pn{2}$. Affine embeddings
always exist: take an arbitrary embedding (with associated closure)
 $Y \subseteq \widehat{Y} \subseteq  \pn{M}$ 
into  some projective space
and blow up the boundary $\widehat{Y} \setminus Y$;
the resulting projective variety ${\overline{Y}}$ contains $Y$
and the complement is a Cartier divisor, so that $Y \subseteq   {\overline{Y}}$
is an affine embedding; finally embed ${\overline{Y}}$ into some projective
space $\pn{N}$: this  embedding is affine. 
If the embedding is not chosen to be   affine,  then the  conclusion of Theorem \ref{swl}  is false, 
as it is illustrated by the example of  the
punctured plane.

We need the following standard vanishing result  due to M. Artin.
\begin{tm}
\label{cdav}
Let $Y$ be an affine variety and $Q \in \pe_Y$ be a perverse sheaf on $Y$.
Then
\[
H^r (Y, Q)  =  0, \; \forall \, r >0, \qquad
H^r_c (Y, Q)  =  0, \; \forall \, r <0.
\]
\end{tm}
{\em Proof.} See \ci{sga4} and \ci{bbd}.\blacksquare

\bigskip
Let $\Lambda \subseteq \pn{N}$ be a  hyperplane, 
 $H := Y \cap \Lambda \subseteq Y$ be the corresponding hyperplane section
 and  consider  the corresponding
open and closed immersions. 
\[ H \stackrel{i}\lorw Y  \stackrel{j}\longleftarrow U: = Y\setminus H. \]
The following is Beilinson's version of the Lefschetz Hyperplane
Theorem.  The proof is an application of  the base change Proposition
\ref{precise01}. 
One can also invoke (as in \ci{be}, Lemma 3.3) the generic base change theorem and reach the same conclusion
(without specifying how one should choose
the hyperplane).
 \begin{tm}
\label{swl}
Let $Q \in \pe_{Y}$ be a perverse sheaf on $Y$.
If $\Lambda$  is a  general hyperplane (for the given affine embedding
$Y \subseteq \pn{N}$), then
$$
H^{r}(Y, j_{!}j^{!}Q) =0, \;\; \forall r < 0, \qquad
H^{r}_c(Y, j_{*}j^{*}Q) =0, \; \forall r > 0.
$$
Moreover, if $Y$ is affine, then 
$$
H^{r}(Y, j_{!}j^{!}Q) =0, \;\; \forall r \neq  0, \qquad
H^{r}_c(Y, j_{*}j^{*}Q) =0, \; \forall r \neq 0.
$$
\end{tm}
{\em Proof.} 
The idea of  proof   is to identify the   cohomology groups
in question
(cohomology groups with compact supports, resp.)
 with
cohomology   groups with compact supports  (cohomology groups , resp.) 
on an auxiliary  {\em affine} variety, and then apply Artin 
vanishing Theorem \ref{cdav}. 

\n Let  $\overline{Y} \subseteq \pn{N}$ be the closure of $Y$.
We have the following chain of equalities (see (\ref{su}))
$$
H^r (Y, j_! j^* Q)= H^r (\overline{Y},J_* j_! j^* Q) \stackrel{=}\longleftarrow
H^r(\overline{Y}, j_! J_* j^* Q)   = H^r_c(\overline{Y}, j_! J_* j^! Q)
=   H^r_c(\overline{U},  J_* j^! Q),
$$
where we have  applied Lemma \ref{pap} and Proposition \ref{precise01}
(applied to $i,J$)
to obtain the second equality.
Since  $j$ is an open immersion, $j^!=j^*$ and
$j^! Q$ is perverse. Since $J$ is an {\em affine} open immersion,
 $J_* j^! Q$ is perverse.
Since $\overline{U}$ is affine, the last group 
is zero for $r >0$ by virtue of Theorem \ref{cdav} and the assertion in cohomology is proved.

\n
The assertion for  $H_c(Y, j_* j^* Q)$,  is proved
 in a similar way.
The relevant sequence of identifications and maps is
$$
H_c(Y, j_* j^* K) = H_c(\overline{Y}, J_! j_* j^* K) \stackrel{=}\lorw
H_c(\overline{Y}, j_* J_! j^* K) = H(\overline{Y}, j_* j^* J_! K)=
H(\overline{U}, j^* J_! K).
$$
\blacksquare

\subsubsection{ A variant of Theorem \ref{swl} using two hyperplane sections}
\label{stuths}
In this section  $Y$ is a  quasi projective variety and we fix an {\em affine}
embedding $Y \subseteq \pn{N}$. 

Let $\Lambda, \Lambda' \subseteq \pn{N}$ be two  hyperplanes, 
 $H := Y \cap \Lambda \subseteq Y$ 
 and  $j: Y \setminus H := U\to Y \leftarrow H : i$  be the
 corresponding
open and closed immersions. 
Note that $j^!=j^*.$ Similarly,  we have $\Lambda', H',  U', i',j'$.
We have the Cartesian diagram of open embeddings 
\begin{equation}
\label{sud}
\xymatrix{
 U    \ar[d]^j & U\cap U'   \ar[l]_{j'}\ar[d]^j  \\
Y &  \ar[l]_{j'}    U'.
}
\end{equation}
Since the embedding 
$Y \subseteq \pn{N}$ is affine, these open embeddings are affine and so are the
open sets $U, U', U\cap U'$. If the embedding were not affine, these
open sets may fail to be affine and the  conclusions 
on vanishing  of Theorem \ref{swl-1}
would not hold.
 Using the  natural maps and isomorphisms (\ref{nmnbc}, \ref{bc1}, \ref{bc2})
and  that $j^! =j^*$, ${j'}^! = {j'}^*$, we get
the following maps
\[ 
j_! j^! j'_* {j'}^* \stackrel{=}\lorw j_! j'_* j^! {j'}^* \stackrel{=}\lorw
j_! j'_* {j'}^* j^!  \stackrel{c'}\lorw j'_* j_! {j'}^* j^! \stackrel{=}\lorw
j'_* {j'}^* j_! j^!
\]
whose composition we denote by
\begin{equation}
\label{wecallc}
c \, : \, 
j_! j^! j'_* {j'}^*  \lorw  j'_* {j'}^* j_! j^!.
\end{equation}
The octahedron axiom applied to the composition
\begin{equation}
\label{eq-200}
j'_* {j'}^* \lorw i_* i^* j'_* {j'}^*  \stackrel{\psi}\lorw i_* j'_* i^* {j'}^*  (= j'_* i_* i^* {j'}^*) 
\end{equation}
yields a distinguished triangle
\begin{equation}
\label{eqn-201}
j_! j^! j'_* {j'}^*  [1] (= j_! j'_* j^! {j'}^*  [1] ) \lorw j'_* j_! j^! {j'}^* [1] \lorw
\mbox{Cone} (\psi) \lorw
\end{equation}
\begin{lm}
\label{pap-1}
The map $$c:j_! j^! j'_* {j'}^* \lorw  j'_* {j'}^* j_! j^!$$ is an isomorphism
if and only if the base change map $$i^* j'_* {j'}^* \lorw j'_* i^* {j'}^*$$ is an 
isomorphism.
\end{lm}
{\em Proof.} Same as Lemma \ref{pap}.
\blacksquare

\begin{tm}
\label{swl-1} Let $Q \in \pe_Y$.
If $(\Lambda, \Lambda')$ is  a general pair, 
then we have 
\[j_! j^! j'_* {j'}^* Q = j'_* {j'}^*   j_! j^! Q \]
and
\[
H^r (Y,  j_! j^! j'_* {j'}^* Q)= H^r_c(Y,  j'_* {j'}^*  j_! j^! Q) = 0, 
\;\; \forall r \neq 0.
\]
\end{tm}
{\em Proof.}
For a fixed  and arbitrary $\Lambda'$,
by virtue of Lemma \ref {pap-1} and Proposition \ref{precise01}
(applied to $i,j'$),
the first equality  holds for
 $\Lambda$ general. This implies that the first equality holds
 for a general pair.
 
 \n
 We prove the statement in cohomology. The one in cohomology with compact supports
 is proved in a similar way, by switching the roles of the two hyperplanes
 $\Lambda$ and $\Lambda'$.
 Note that $j'_* {j'}^* Q$ is perverse. The vanishing of the groups for $r <0$
 and $\Lambda$ general
 follows from Theorem \ref{swl}. The vanishing for $r>0$
 is obtained as follows: $H^r (Y,  j_! j^! j'_* {j'}^* Q) = 
 H^r (Y,  j'_* {j'}^* j_! j^!   Q) =  H^r (U',   {j'}^* j_! j^!   Q)$; $U'$ is affine,
 ${j'}^* j_! j^!   Q$ is perverse
 and the last  group is zero for $r <0$ by Theorem \ref{cdav}.
 \blacksquare
 
 \begin{rmk}
 \label{usato}
 {\rm Theorem \ref{swl-1} is due to Beilinson \ci{be} and it is used in
 \ci{decmigs1} to describe perverse filtrations on quasi projective varieties
 using general pairs of flags.}
 \end{rmk}

\subsubsection{The Lefschetz hyperplane theorem for constructible sheaves}
\label{fshe}
As it is observed in \ci{nori}, Introduction, Theorem \ref{swl} admits a sheaf-theoretic version
which we state and prove below.
Let $Y \subseteq \pn{N} $ be a quasi projective  variety of dimension $n$ embedded in some projective space in such a way 
 that the embedding is affine. Let ${\cal V} \subseteq \pn{N}$ be a hypersurface,
$V: = Y \cap \cal V$ and $j: Y \setminus V \to Y$.

\begin{tm}\label{swls}
Let $T$ be a constructible sheaf on $Y$. There is a hypersurface
${\cal V}$ such that
\begin{enumerate}
\item
$H^r(Y, j_! j^! T) =0$, for every $r < n$ (for every $r \neq n$ if $Y$ is affine),
\item
$\dim{V} = \dim{Y} -1$.
\end{enumerate}
\end{tm}
{\em Proof.}
Let $\Sigma$ be a stratification of $Y$ with respect to which
$T$ is constructible.   The union $S_n$
of all $n$-dimensional strata is a non-empty, Zariski open
subvariety of $Y$ with the property that $T_{|S_n}$ is locally constant
and $\dim{(Y \setminus S_n)} \leq n-1$.  Let ${\cal V}'\subseteq \pn{N}$ be a hypersurface
containing $Y \setminus S_n$ but not containing any of
the irreducible components of $S_n$. Since
the open embedding   $j': Y\setminus V' \to Y$
is affine, $j'_! {j'}^! T[n]$ is a perverse  sheaf on $Y$.
We apply Theorem \ref{swl} to this perverse sheaf  and conclude
that the desired hypersurface is of the form
${\cal V}:= {\cal V}' \cup \Lambda$ for some general
hyperplane $\Lambda$. \blacksquare

\begin{rmk}
\label{hh}
{\rm
The hypersurface ${\cal V}'$ must contain the ``bad locus" of the sheaf $T$. In particular,
it is a ``special" hypersurface of sufficiently  high degree.  As the proof shows, it is not necessary 
to achieve 2. in order to achieve 1. However, 2. is useful in procedures
where one uses induction on the dimension. I do not know of a version
of Theorem \ref{swls}
for cohomology groups with  compact supports.}
\end{rmk}

\subsection{The generic base change theorem}
\label{s-gbc}
The Generic Base Change Theorem was proved in \ci{sga4m} as an essential ingredient,
in the \'etale context,   towards the constructibility for direct images
of complexes with constructible cohomology sheaves  for
morphisms of finite type over a field.
These kinds of constructibility results are fundamental and permeate the whole
theory of 
\'etale cohomology.

In this section we want to state the Generic Base Change Theorem
and show how it can be applied in practice when one has a base change
issues with ``parameters," e.g. elements of a linear system. 
For example, in the proof of the Lefschetz hyperplane  theorem \ref{swl}
if one can afford to work with general linear sections, then 
Proposition \ref{pr-gbcf} can be used in place of  Proposition \ref{precise01}.

\subsubsection{Statement of the generic base change theorem}
\label{sgbc}
Let $X \stackrel{f}\to Y \stackrel{p}\to S$ and  $S' \to S$ be  maps.
Denote by  $X' \stackrel{f}\to Y' \stackrel{p}\to S'$ the
varieties and maps obtained by base change via the given $S' \to S.$

Let $C \in {\cal D}_X.$ One says that {\em the formation of $f_*C$
commutes with arbitrary base change}
if, for every $S' \to S,$ the  resulting first base change map (\ref{bc1})
$g^* f_* C \to  f_* g^* C$  
is an isomorphism. 

The issue  does not arise  for $f_!$, in fact,
 the base change isomorphism (\ref{bc2})
$g^*f_! = f_!g^*$
implies that for every $C \in {\cal D}_X$ the formation
of $f_!C$  commutes with arbitrary  base change.

Given a stratification $\Sigma$ of $X$ and a map $f: X \to Y,$
one can refine $\Sigma$ so that the refinement is part of a stratification
of the map $f.$ 
It follows at once that, given $f: X \to Y$ and $C \in {\cal D}_X,$
there is an Zariski-dense open subset
$U \subseteq Y$ with the property that, given $f^{-1} (U) \to U =U,$
the formation of $f_* (C_{ | f^{-1}U} )$
commutes with arbitrary base change.
It is sufficient to take for $U$ the dense open stratum on $Y$
of the stratification for $f$ refining the one for $C$. In fact, $f$ is then 
topologically locally trivial over $U$ and the base change maps are then isomorphisms.

However, what above is insufficient to prove the, for example,
the   vanishing Theorem
\ref{swl}. Moreover, it cannot be used for example to work with
constructible sheaves for the 
\'etale topology for varieties over a field, where one cannot achieve the local triviality
of $f: X \to Y$ over $U \subseteq Y$ (in fact, the generic base change theorem
is a tool that effectively fixes this problem at the level of sheaves).

The Generic Base Change Theorem is a tool apt to deal with these and other situations.

\begin{tm}
\label{gbc} 
Let $X \stackrel{f}\to Y \stackrel{p}\to S$ be maps and $C \in {\cal D}_X.$ There 
exists a Zariski open and dense  subset $V \subseteq S$ such that,
if one takes $ (pf)^{-1} (V) \to  p^{-1} V \to V,$ then 
the formation of $f_*(C_{ | (pf)^{-1} V } ) $ 
commutes with arbitrary base change $T \to V.$ 
\end{tm}
{\em Proof.} For  the \'etale  case see \ci{sga4m}, [Th. finitude], Th. 1.9.
The proof in the case  of complex varieties and $C \in {\cal D}_X$ is similar
and, in fact, simpler.
\blacksquare

\begin{rmk}
\label{sigmank}
{\rm
Note that the open set $V$ depends on $C$. However, an inspection of the proof
reveals that given a stratification $\Sigma $ of $X$, one can choose
the Zariski open and dense subset $V \subseteq S$ so that the conclusion
of the Generic Base Change Theorem holds for every $C' \in {\cal D}_X^{\Sigma}$.
}\end{rmk}
\begin{rmk}
\label{rmk-33}
{\rm 
For $\Sigma$ and $V$ as above, the formation of $f_!$ commutes with the formation
of $g^!$ over $V$ for every $K \in {\cal D}_X^{\Sigma}$.
In fact, to prove that $g^!f_! K  \leftarrow f_! g^! K$ is an isomorphism
for $K \in D^{\Sigma}_X$,
it is sufficient to observe that $K^{\vee} \in D^{\Sigma}_X$
and  dualize the isomorphism $g^*f_* K^{\vee}  \to
 f_* g^* K^\vee$ which holds over $V$ by Theorem \ref{gbc}.
}
\end{rmk}

\subsubsection{Generic base change theorem and families of hyperplane sections}
\label{gbcfh}
The following standard  lemma is an illustration of the use of Generic Base Change.
It is essentially a special case, formulated in a way that
directly applies to the situation dealt-with in the Lefschetz Hyperplane
Theorem.

\begin{lm}
\label{bct}
Let $f: X \to Y$ be a map, 
$C \in {\cal D}_X$ and 
$$
\xymatrix{
X_{2,T}  \ar[r]^{\tau'}  \ar[d]^{f_{2,T}}  
\ar@/^2pc/[rrrr]^{i'_T}
&  {X}_{2,V}\ar[d]^{f_{2,V}}   \ar[r]^{u'} 
&  {X_2} \ar[d]^{ f_2 } \ar[r]^{v'}  
\ar@/^2pc/[rr]^{g'}
 &  {X_1}  \ar[d]^{ f_{1} } \ar[r]^{\pi'} 
& X \ar[d]^{f}
 \\
{Y_{2,T}}  \ar[dd]^{p_T}   \ar[r]^{{\tau}}  
 \ar@/_2pc/[rrrr]^{{i}_T}
 &  { Y_{2,V}}  \ar[dd]^{p_V} 
 \ar[r]^{ {u} }  &Y_2 \ar[dd]^p  \ar[r]^{{v} }  
 \ar@/_2pc/[rr]_{{g}}
 &  Y_1 \ar[r]^{ {\pi} } 
 & Y 
 \\
 &&&&\\
 T  \ar[r]^t & V \ar[r]^{u''}  & {S}
 }
$$
be a commutative diagram with Cartesian squares satisfying:
\begin{enumerate}
\item
${\pi}$ smooth; in particular,  $ {\pi}^* f_* \simeq {f_{1}}_* {\pi'}^*$;
\item
${g}$ smooth; in particular, ${g}^* f_* \simeq {f_2}_*  {g'}^*$;
\item 
$V \subseteq S$ is a Zariski-dense open subset such that the formation 
of ${f_{2,V}}_* {u'}^*{v'}^* {\pi'}^*C$ 
commutes with arbitrary base change (on $V$).
\end{enumerate}
For every $t: T \to V,$  the natural  base change map
 is an  isomorphism:
$$
{i}_T^*  f_* C  \stackrel{\simeq}\lorw {f_{2,T}}_* {i'}^*_T C.
$$
\end{lm}
{\em Proof.}
The natural map
$$
{i}_T^*  f_* C  \lorw {f_{2,T}}_* i^*_T C
$$
factors as follows
$$
{i}_T^*  f_* C =  {\tau}^*{u}^*{g}^*f_*C \lorw
{\tau}^* {u}^*  {f_{2}}_*  {g'}^* C  \lorw
{\tau}^*  {f_{2,V}}_* {u'}^* {g'}^* C \lorw
{f_{2,T}}_* {\tau'}^* {u'}^* {g'}^* C =
{f_{2,T}}_* {i'}^*_T  C.
$$
Since ${g}$  
and ${u}$ are smooth, the first  and second arrows are isomorphisms
The third one is an isomorphism by the choice of $V.$
\blacksquare

\begin{rmk}
\label{rmk-1}
{\rm
Fix a stratification $\Sigma$ for $X$. As in Remark \ref{sigmank}, 
we can choose $V$ so that  3. above holds for every $C' \in {\cal D}_X^{\Sigma}$
and conclude (see Remark \ref{rmk-33}) that we have the base change isomorphisms
\[
{i}_T^*  f_* C'  \stackrel{\simeq}\lorw {f_{2,T}}_* {i'}^*_T C',
\qquad
{i}_T^!  f_! \,C'  \stackrel{\simeq}\longleftarrow{f_{2,T}}_! {i'_T}\!\!^! \,C',
\qquad \forall \, T \lorw V, \;\; \forall\, C' \in D^{\Sigma}_Y.
\]
}
\end{rmk}

We now  apply Lemma \ref{bct} and Remark \ref{rmk-1} to the following
situation:
let $f: X \to Y$ be a map of varieties, $|{\frak H}|$  be a finite dimensional 
and base-point-free linear system on $Y$, e.g. the very ample linear system
associated with an embedding on $Y$ into projective space.
Given $H \in |{\frak H}|$, we have  the Cartesian diagram
(\ref{cdm-000}).
\begin{pr}
\label{pr-gbcf}
Let $\Sigma$ be a  stratification of $X$.
If $H \in |{\frak H}|$ is  general,  then for every $C \in D^{\Sigma}_X$
 the base change maps 
\[ i^* f_* C \lorw f_* i^* C ,\qquad   i^! f_! C\longleftarrow f_! i^!  C\]
are isomorphisms.
\end{pr}
{\em 1st proof } (it uses the generic base change theorem
and it does not single out a specific $H$). We only need to apply Lemma \ref{bct} and Remark \ref{rmk-1}
to the following situation: $Y_1 := Y \times |{\frak H}|$, $Y_2 \subseteq Y \times
|{\frak H}|$ the universal hyperplane section, $S:= |{\frak H}|$ and 
$t: T \to V$ is the  embedding of a closed  point.

\smallskip
\n
{\em 2nd proof}  (it uses Proposition  
\ref{precise01} and it identifies precisely which conditions on $H$
must be satisfied). A general hyperplane $\Lambda$ is transverse
to all the strata of a fixed stratification  of  $f_*C$. This means
that $i:H \to Y$ is a normally nonsingular inclusion
with respect to the given stratification. For such a $\Lambda$,
$i^* f_* C = i^![2] f_*C= f_* i^! C [2]$. For $\Lambda$ general,
$i: X_H \to X$ is transverse to all the strata of a fixed stratification
of $C$ and we have $i^! C = i^* C [-2]$. This establishes that the first base change map
in question is an isomorphism. The proof for the second one is
similar.
\blacksquare

\end{document}